\documentclass{llncs} 

\usepackage{amssymb,amsfonts}
\usepackage{xcolor}
\usepackage{soul}

\usepackage{cite}
\usepackage{url}
\usepackage{bm}
\definecolor{ggreen}{rgb}{0,.47,0}
\newcommand{\TheoCite}[2]{\par\vspace{.95em}\noindent\textbf{#1.\ }{\it #2}\vspace{.95em}}    %

\def\xx#1{} 
\def\xy{\hspace{.11em}}                                         
\def\xxy{\hspace{.14em}}                                         
\def\xz{\hspace{-.11em}}                                        
\def\xxz{\hspace{-.05em}}                                        
\def\a{{\alpha}}
\def\b{{\beta}}
\def\x{\bm x}  
\def\y{\bm y}
\def\z{\bm z}
\def\p{\bm p}
\def\q{\p'}
\def\ms{\mathstrut}
\def\msd{_{\mathstrut}} 
\def\msu{^{\mathstrut}}

\def\aaa#1{f\xxz_{\bm#1\msd}^{P}\xz}
\def\bbb#1{f^{\bm#1\ms}_{\!P}}
\def\aa{\aaa{\uparrow}}
\def\bb{\bbb{\downarrow}}
\def\ada{\aaa{\upharpoonright}}
\def\bdb{\bbb{\downharpoonright}}
\def\aal{f\xxz_{\bm\uparrow}^{P}\!\!\xz} 
\def\bbl{f^{\bm\downarrow}_{\!P}}       
\def\aau{f\xxz_{\bm\uparrow\msd}^{P\msu}\xz}

\def\ax{\aa(\x)}
\def\bx{\bb(\x)}
\def\adx{\ada(\x)}
\def\bdx{\bdb(\x)}
\def\axl{\aal(\x)}
\def\axu{\aau(\x)}
\def\bxl{\bbl(\x)}

\def\ap{\aa(\p)}
\def\bp{\bb(\p)}

\def\u01{u\zaz{01}\xz}
\def\uab{u\zaz{\a\b}}
\def\uabx{\uab(\x)}
\def\fx{f(\x)}
\def\fo{f}
\def\fhat{\breve{f}} 

\def\zaz#1{\mathop{\!^{}_{#1}}}

\def\cdc{,\ldots,}
\def\ok{1\cdc k}
\def\fP{f\zaz{\!\!P}}
\def\fX{f\zaz{\!\!X}}
\def\fPi{f^{-1}_{\!P}}
\def\fPx{\fP(\x)}
\def\fPp{\fP(\p)}
\def\fPq{\fP(\q)}

\def\fPo{\fP}
\def\X{X}
\def\Xm{X_{\bm\uparrow}^{}\!\!\xz}
\def\Xp{X^{\bm\downarrow}}
\def\Pm{P_{\bm\uparrow}^{}\!\!\xz}   
\def\Pp{P^{\bm\downarrow}} 
\def\Pmx{\Pm(\x)}
\def\Ppx{\Pp(\x)}
\def\Pmp{\Pm(\p)}
\def\Ppp{\Pp(\p)}
\def\Pmr{\fPi(\xxz(-\infty, r])} 
\def\Ppr#1{\fPi(\xxz[\xy#1,+\infty))} 
\def\precE{\preccurlyeq}
\def\succE{\succcurlyeq}
\def\succEe{\succE} 
\def\succN{\succ} 

\def\succEP{\succE\!\zaz{P}}
\newcommand{\eq}[2]{\begin{equation}\label{#1}#2\end{equation}}

\def\To{\Rightarrow}
\def\ToTo{\Leftrightarrow}

\def\painfty{\bm{+\infty}}
\def\mainfty{\bm{-\infty}}

\def\Xr{{(\X,\succE)}}
\def\XX{(\X,\tau,\succE)}

\def\R{{\mathbb R}}
\def\N{{\mathbb N}}
\def\Z{{\mathbb Z}}
\def\RR{Y}
\def\Rk{\R^k}
\def\tR{\widetilde{\R}}
\def\tX{\widetilde{X}}

\def\uppercontour{upper $P$-contour} 
\def\lowercontour{lower $P$-contour} 
\def\Vx#1{\mathcal{V}_{#1}^{\x}}
\def\eqeqA{\big[\xy\p,\,\q\in P\mbox{ and\, }\q\approx\p\xy\big]\,\To\fPq=\fPp}

\DeclareFontFamily{U}{mathx}{}
\DeclareFontShape{U}{mathx}{m}{n}{ <-> mathx10 }{}
\DeclareSymbolFont{mathx}{U}{mathx}{m}{n}
\DeclareFontSubstitution{U}{mathx}{m}{n}

\DeclareMathAccent{\widecheck}{0}{mathx}{"71}

\newcommand{\keywords}{\vspace*{1ex}{\bf Keywords: }}

\begin{document}
\mainmatter
\title{Extending Utility Functions on Arbitrary Sets
}
\titlerunning{Extending Utility Functions on Arbitrary Sets\xx{ without Continuity Constraints}}
\date{April 8, 2001}

\author{Pavel Chebotarev\xx{\thanksref{coord}}}
\institute{
Technion--Israel Institute of Technology, Haifa, 3200003 Israel;\\ 
\xx{and }A.A.~Kharkevich Institute for Information Transmission Problems, RAS, 19 Bol'shoi Karetnyi per., Moscow, 127051 Russia,
{\em E-mail\/}: pavel4e@technion.ac.il \xx{gmail.com}
}
\authorrunning{Pavel Chebotarev}
\maketitle

\begin{abstract}
We consider the problem of extending a function $\fP$ defined on a\xx{n arbitrary} subset $P$ of an arbitrary\xx{ preordered} set~$\X$ to~$\X$ strictly monotonically with respect to a preorder $\succE$ defined on $\X$, without imposing continuity constraints. 
We show\xx{prove} that whenever $\succE$ has a utility representation, $\fP$ is extendable if and only if it is gap-safe increasing.
A class of extensions involving an arbitrary utility representation of $\succE$ is proposed and investigated\xx{studied}.
Connections to related topological results are discussed.
The\xx{ necessary and sufficient} condition of extendability and the form of the extension\xx{proposed construction} are simplified when
$P$ is a Pareto set\xx{subset of~$X$}.

\keywords
Extension of a utility function; Monotonicity; Utility representation of a preorder; Pareto set
\end{abstract}


\section{Introduction}
\label{s:Intro}

Suppose that a decision maker has a utility function $\fP$ defined on some subset $P$ of a Euclidean space $\Rk$ of alternatives.
It is usually assumed that $\fP$ strictly increases in $k$ coordinates corresponding to particular criteria. Therefore, it is of interest\xx{a natural problem} to determine\xx{  establish find} conditions under which\xx{ a strictly increasing (this) function defined on $P$} this function can be\xx{ strictly monotonically (preserving increase)} extended in such a way as to provide a strictly increasing function on~$\Rk$. In settings\xx{situations} where all elements of $\Rk$ are theoretically feasible, these conditions can be considered as those of the consistency of~$\fP$.

In this paper, we consider a\xx{the} more general problem where an arbitrary preordered set $\Xr$ is substituted for~$\Rk$. We show that whenever preorder $\succE$\xx{ has} enables a utility representation (in the form of a real-valued function strictly increasing w.r.t. $\succE$ on~$\X$), a strictly increasing extension of $\fP$ to $\Xr$ exists if and only if $\fP$\xx{ this function} is gap-safe increasing with respect to~$\succE$.
Moreover\xx{Furthermore}, such an extension can be based on any utility representation of~$\succE$. The main object of this paper is the general class of extensions (\ref{f}) involving an arbitrary $(0,1)$-utility representation $\u01$ of $\succE$ and arbitrary real constants $\a$ and $\b>\a.$

\xx{Furthermore1/20}We also consider the\xx{ special} case where the structure of subset $P$ restricts\xx{constrains1/2 the} functions strictly increasing on $P$ to the minimum extent. This is the case where $P$ is a Pareto set; for such $P,$ the extension takes\xx{ on} a simpler form.
\xx{Corollary~\ref{c:Par} addresses this situation.}

Starting with the classical\xx{ papers} results of Eilenberg \cite{Eilenberg41OrderTopo}, Nachbin \cite{Nachbin50HahnBa,Nachbin65TopolOrderNo}, and Debreu \cite{Debreu54Repr,Debreu59Topol,Debreu64Conti}, much of the work related to utility functions has been done under the continuity assumption~\xx{ that they are continuous}\cite{BridgesMehta95ReprPrefOrd,Evren21Extension}. \xx{Often(Sometimes)}In some cases this assumption is made ``for purposes of mathematical reasoning''~\cite{Allen34Indiffer}.
On the other hand\xx{However}, this requirement is not always necessary. \xx{Furthermore}Moreover, there are threshold effects \cite{Hirshleifer92Uncertainty,Kritzman94TimeDivers} such as a\xx{ transition} shift from quantity to quality or disaster avoidance behavior that require utility jumps. In other situations, the feasible set of possible outcomes is a discrete or finite subset of the entire\xx{ whole} space, which may eliminate or relax the continuity\xx{ requirement} constraints.
\xx{As a result Therefore}Thus, utility functions that\xx{ are not everywhere continuous} may not be continuous everywhere are useful or even necessary to model some\xx{a number of} real-world\xx{applications} problems~\cite{GaleSutherland68Good,Masson74Discontin,Hande07Inelastic,Diecidue08Aspirat,Siciliani09paying,Andreoni10Allais,Bian19Discontin}. For a discussion of various versions of the continuity postulate in utility theory, we refer to~\cite{Uyanik22continuity}.

In this paper, we study\xx{ consider focus on} the problem of extending utility functions defined on arbitrary subsets of an arbitrary set $X$ equipped with a preorder $\succE$, but not endowed with a topological structure, since we do not impose continuity requirements.
However, some kind of continuity of an associated\xx{related} inverse mapping\xx{-ping appears (acts) as a} follows from the necessary and sufficient condition of\xx{ its} extendability we\xx{ obtain present} establish.

\section{The problem and standard definitions} 
\label{s:ProbDef}

Throughout the paper $\Xr$ is a {\em preordered set\/}, where $\X$ is an arbitrary nonempty set and $\succE$ is a {\em preorder\/} (i.e., a transitive and reflexive binary relation) defined on~$\X.$
We first formulate the problem under consideration and then provide the necessary definitions; the definitions of basic\xx{relevant} properties and classes of binary relations are given in Appendix~\ref{a:BinaryR}\xx{; for convenience, some of them are mentioned in the main text}.

Consider any subset $P\subset \X$\xx{ of $\xy\Rk$} and any\xx{ strictly increasing (with respect to $\succEe$)} real-valued function $\fPo$
defined on~$P$. The\xx{ main} problem studied in this paper is:
to find\xx{ a necessary and sufficient} conditions under which $\fPo$ can be\xx{ strictly monotonically} extended to $\Xr$ yielding\xx{producing} a strictly increasing function and to construct a fairly general class of such extensions when they exist.

\xx{We now recall}
The definitions of the\xx{some (the main)} relevant terms are as follows.

Given a preorder $\succE$ on $X,$ the {\em asymmetric $\succ$} and {\em symmetric $\approx$ parts\/} of $\succE$ are the relations\footnote{The elements of $X$ are printed\xx{shown} in bold (as is common for vectors in~$\Rk$) to distinguish them from real numbers.}
$[\x\succ  \y]\equiv    [\x\succE \y$ and not $\y\succE\x$] and
$[\x\approx\y]\equiv    [\x\succE \y$ and     $\y\succE\x$], respectively,
where $\equiv$ means ``identity by definition.''
Relation $\succ$ is transitive and irreflexive (i.e., it is a {strict partial order\xx{\footnote{A counterpart of a strict partial order is a {\em partial order}, i.e., a transitive, reflexive and antisymmetric relation, or, equivalently\xx{in other words}, an antisymmetric preorder.}}}), whereas $\approx$ is transitive, reflexive, and symmetric (i.e., an {equivalence relation}).

The \emph{converse relations\/} corresponding to $\succE$ and $\succ$ are
$\precE\,:[\x\precE\y]\equiv[\y\succE\x]$ and
$ \prec\,:[\x\prec \y]\equiv[\y\succ \x]$.

For any $P\subseteq X,$ $\succE\!\zaz{P}$ is the restriction of $\succE$ to $P$.

%
%
$\x\in X$ is a {\em maximal\/} (\xz{\em minimal\/}) {\em element\/} of $\Xr$ iff $\x'\succ\x$ (resp., $\x\succ\x'$) for no $\x'\in X.$

\begin{definition}
\label{d:defi}
{\rm
A function $\fP\!:P\to\R,$ where $P\subseteq X,$ is said to be {\em weakly increasing with respect to the preorder $\succEe$} defined on $\X$ (or, briefly, {\em weakly increasing\/}) if for all $\p,\,\q\in P,\,$ $\q\succE\p$ implies\footnote{In a different terminology \cite{Birkhoff40Lattice}, functions with this property are referred to\xx{ called} as {\em order-preserving\/}, or {\em isotone.}}
${\fPq\ge\fPp.}$

If, in addition, $\fPq>\fPp$ for all $\p,\,\q\in P$ such that $\q\succ\p,$ then $\fP$ 
is called {\em strictly increasing\/} w.r.t. $\succEe$, or a {\em utility representation of\/}~$\succEP$.
}
\end{definition}

Utility functions $\fPo$ strictly increasing w.r.t. $\succEe$ 
can express the attitude, consistent with the preference preorder $\succEe$, of a decision maker towards the elements of~$P$.
Utility representations of preorders and partial orders have been studied since~\cite{Nachbin65TopolOrderNo,Aumann62Utility,Peleg70Utility,Fishburn70Utility}.

\smallskip
It follows from Definition~\ref{d:defi} that for any weakly increasing function $\fP$,
\eq{e:appro}{
\eqeqA
}
Using (\ref{e:appro}) we obtain the following simple lemma.

\begin{lemma}
\label{l:strict}
A function $\fP\!:P\to\R,$ where $P\subseteq X,$ is strictly increasing with respect to a preorder $\succEe$ defined on $\X$ if and only if for all\/ $\p,\,\q\in P,\,$
\eq{e:l:strict}{
\big[\xy\q\approx\p\:\To\fPq=\fPp\big] \,\mbox{ and }\;\big[\xy\q\succ\p\:\To \fPq>\fPp\big],
}
where $\approx$ and $\succ$ are the symmetric and asymmetric parts of $\succE,$ respectively.
\end{lemma}

Indeed, (\ref{e:l:strict}) follows from Definition~\ref{d:defi} using~(\ref{e:appro}). Conversely, if (\ref{e:l:strict}) holds, then $[\xy\q\succE\p\:\To\fPq\ge\fPp],$ since $\q\succE\p$ implies $[\xy\q\approx\p\mbox{ or }\q\succ\p]$ with the desired conclusion in either\xx{each} case, while the second condition is immediate.

\begin{definition}  
\label{def2}
{\rm
A real-valued function $\fPo$ defined on $P\subseteq X$ is {\em strictly monotonically\/\footnote{We mean increase.} extendable\xx{ ``extendible" is also allowed} to\/ $\Xr$} if there exists a function $\fo\equiv\fX:\;\X\to\R$ such that
\\$(\ast)$ the restriction of $\fo$ to $P$ coincides with $\fPo,$ and
\\$(\ast\ast)$ $\fo$ is strictly increasing on $\X$ with respect to~$\succEe$.\\
In this case$,$ $\fo$ is said to be a {\em strictly increasing extension\/} of $\fPo$ to~$\Xr$.
}
\end{definition}

In economics and decision making, \xx{objects/}alternatives are often\xx{ characterized by specific criteria and} identified with $k$-dimensional vectors of criteria values~\cite{Thakkar21MCDM} or\xx{ vectors of} goods~\cite{Allen34Indiffer}. In such cases, $X=\Rk.$ Thus, an important special case of the extendability\xx{ extension} problem is the problem of extending to~$\Rk$ functions defined on $P\subset\Rk$ and strictly increasing w.r.t. the Pareto preorder on~$\Rk$.
The {\em Pareto preorder\/} $\succE$ \cite{Debreu93SmalePareto} is defined as follows:
for any $\x=(x_1\cdc x_k)$ and $\y=(y_1\cdc y_k)$ that belong to $\Rk,$ $[\x\succE \y]\equiv [x\zaz{i}\ge y\zaz{i}$ for all $i\in\{\ok\}]$.

\section{Extensions of preorders and corresponding utilities}
\label{s:pre&util}

Extensions of preorders and partial orders and their numerical representations have been studied since\xx{ Zorn's lemma and} Szpilrajn's theorem \cite{Szpilrajn30Exten} according to which every\xx{ strict} partial order can be extended\xx{extends} to a\xx{ strict} linear order.

Another basic result is that a preorder $\succE$ has a utility representation whenever\xx{ provided that (*when) under the condition\xx{s such as} of the existence of} there exists a countable dense\footnote{$Y\subseteq\X$ is $R$-{\em dense\/} in $\X,$ where $R$ is a binary relation on $\X$ \cite{Fishburn70Utility}, iff $\x'R\x$ $\To$ [$\x'\in Y$ or $\x\in Y$ or [$\x'R\y$ and $\y R\x$ for some $\y\in Y$]] for all $\x,\x'\in\X$.} (w.r.t.\ the induced partial order) subset in the factor set $X/\!\approx,$ where $\approx$ is the symmetric part of $\succE$ \cite{Debreu64Conti,Richter66,Fishburn70Utility}. This\xx{It} is not a necessary condition\xx{ of the representability}, however, for the subclass of {weak orders\/} (i.e., connected preorders\xx{transitive and complete relations}), it is necessary.

Among the extensions of the Pareto preorder on $\Rk$ are all lexicographic linear orders \cite{Fishburn70Utility} on~$\Rk$. When $k>1$, these extensions lack\xx{have no} utility representations \cite{Debreu64Conti}, while a utility representation of the Pareto preorder\xx{$\succE$ \cite{Debreu93SmalePareto}} is any function strictly increasing in all coordinates.

\xx{In the general case, all}
Any\xx{ possible} utility representation of a preorder $\succE$\xx{ (if they exist)} induces a\xx{ strict} weak order that extends~$\succE$. In turn, this\xx{such a (strict)} weak order determines its utility representation up to an arbitrary strictly increasing transformation; for certain related results, see~\cite{Fishburn70Utility,Morkeliunas86b,Morkeliunas86a,Herden89Existence,BridgesMehta95ReprPrefOrd,Evren21Extension}.
As was seen on the example of the Pareto preorder\xx{relation}, not all\xx{ strict} weak orders extending $\succE$ correspond to utility representations of~$\succE$. However, this is true when $\X$ is a vector space and the weak order has\xx{ A sufficient condition for this\xx{such a representability (this connection)} is} the Archimedean property, which ensures \cite{Fishburn70Utility} the existence of a countable dense (w.r.t.\ this weak order) subset of~$\X$.

\section{Preliminaries}
\label{s:prelim}

Theorem~\ref{T} below gives\xx{provides} a necessary and sufficient condition for the strictly increasing extendability of a function defined on a subset of $\X$ w.r.t.\ a preorder\xx{ $\succE$ on $\X$} that has a utility representation. Moreover, this theorem \xx{presents}introduces a class of\xx{ such} extensions that depend on both the initial function and an arbitrary utility representation of the preorder.
%
%

\smallskip

We now\xx{Let us} introduce the\xx{ basic} notation used in Theorem~\ref{T} and\xx{ discuss} present simple facts related to it.

Let $\tR$ be the extended real line:
\eq{e:extR}{
\tR=\R\xy\cup\xy\{-\infty,+\infty\}
}
with the ordinary $>$ relation supplemented by $+\infty>-\infty$ and $+\infty>x>-\infty$ for all $x\in\R$.
Since the extended $>$ relation is a strict linear order, it determines unique smallest ($\min Q$) and largest ($\max Q$) elements in any nonempty finite $Q\subset\tR$.

Functions $\sup Q$ and $\inf Q$\xx{will be} are considered as maps from $2^{\R}$ to $\tR$ defined for $Q=\emptyset$ as follows:
$\sup\emptyset=-\infty$ and
$\inf\emptyset=+\infty$.
This\xx{ definition} preserves \emph{inclusion monotonicity\/}, i.e., the property that $\sup Q$ does not decrease and $\inf Q$ does not increase with the expansion of the set $Q$ (cf.~\cite[Section~4]{Tanino88supremum}).
Throughout we assume $+\infty+x=+\infty$ and $-\infty-x=-\infty$ whenever $x>-\infty,$ while indeterminacies like $+\infty+(-\infty)$\xx{ will} never \xx{appear}occur in our formulas.

\begin{remark}\label{r:sup}
If $\RR\subset\R$ and $\RR$ is bounded, then defining $\sup Q$ and $\inf Q$ on\xx{ as functions from} $2^{\RR}$ with the preservation of inclusion monotonicity allows setting $\sup\emptyset=a$ and $\inf\emptyset=b,$ where $a$ and $b$ are any strict lower and upper bounds of $\RR,$ respectively. This\xx{ observation can be used (can be done in)} is applicable to (\ref{axbx}) below\xx{ whenever}\xx{ in the case where} whenever the\xx{ set of values} range of $\fP$ is bounded.\xx{!!Check!}
\end{remark}

\begin{definition}
\label{def23}
{\rm
\xx{Given}For any $P\subseteq\X$ and $\x\in\X,$
the {\em lower $P$-contour\xx{ of $\x$}} and
the {\em upper $P$-contour of\/ $\x$} are\xx{ the sets}
$\Pmx\equiv\{\p\in P\mid\p\precE\x\}$ and
$\Ppx\equiv\{\p\in P\mid\p\succE\x\},$
respectively.
}
\end{definition}

For any $\fPo\!:P\to\R$, where $P\subseteq \X$, define two\xx{ auxiliary} functions from $\X$ to $\tR$:
\eq{axbx}{
\begin{array}{rrl}
\ax=&\sup\big\{\msd \fP(\p)&\mid\xy\p\in \Pmx\big\};\\
\bx=&\inf\big\{\msu\fP(\p)&\mid\xy\p\in \Ppx\big\}.
\end{array}
}
By definition, the ``lower supremum'' $\ax$ and ``upper infimum'' $\bx$ functions\xx{ corresponding to $\fPx$} can take values $-\infty$ and $+\infty$ along with real values.
%

\smallskip
It follows from the transitivity of $\succEe$\xx{ relation} and the inclusion monotonicity of the $\sup$ and $\inf$ functions that for any (not necessarily increasing)\xx{ function} $\fPo,$ functions $\ax$ and $\bx$ are {weakly increasing\/}\xx{} with respect to~$\succEe$:
\eq{e:abmono}{ 
\mbox{For all } \x,\x'\in \X,\,\x'\succEe\x \mbox{ implies } \big[\aa(\x')\ge\ax\mbox{ and }\bb(\x')\ge\bx\big].
}

Consequently,
\eq{e:ApprCont}{ 
\mbox{for all } \x,\x'\in \X,\,\x'\approx\x \mbox{ implies } \big[\aa(\x')=\ax\mbox{ and }\bb(\x')=\bx\big].
}

Furthermore,\xx{ Besides,} since $\p\in P$ implies $\p\in\Pmp\cap\Ppp,$
it holds that
\eq{a>f>b}{ 
\mbox{for all }\,\p\in P,\quad\ap\ge\fPp\ge\bp.
}

We will use the following characterizations of the class of weakly increasing functions $\fPo$ in terms of $\aa$ and~$\bb.$

\begin{proposition}
\label{P:nonstr}
For any $P\subseteq\X$ and $\fP\!:P\to\R,$ the following statements\xx{ $(i)$ to $(v)$} are equivalent\/$:$

\noindent
\xxy$(i)$ \,\,$\fPo$ is weakly increasing\/$;$

%
%
%
%
%
\noindent
$\begin{array}{lrll}
(ii) &      \bx&\ge\ax\,  &\mbox{for all }\;\x\in\X;\\
(iii)& \bb(\x')&\ge\ax    \msu&\mbox{for all }\;\x,\,\x'\in\X\mbox{ such that }\x'\succEe\x;\\
(iv) &\!\fP(\p)&\ge\ap    \msu&\mbox{for all }\;\p\in P;\\
(v)  &      \bp&\ge\fP(\p)\msu&\mbox{for all }\;\p\in P;\\
(vi) &      \bp&\ge\ap    \msu&\mbox{for all }\;\p\in P.
\end{array}
$
\end{proposition}

\xx{All}The proofs are given in Section~\ref{a:Proofs}. 

\begin{remark}\label{r:=}
In view of Eq.~(\ref{a>f>b}), inequality in items $(iv)$ to $(vi)$ of Proposition~\ref{P:nonstr} can be replaced by equality.
\end{remark}

\section{Gap-safe increasing functions}
\label{s:sepinc}

In this section,\xx{ for any $P\subseteq\X,$} we\xx{We now} consider the class of gap-safe increasing functions $\fPo,$ which is not wider, but can be narrower for some $\X$ and $P$ than the class of strictly increasing functions $P\to\R$ (see Proposition~\ref{P0} and Example~\ref{ex:f} below). We will show that this is precisely\xx{exactly} the class of functions that admit\xx{allow, enable} a strictly increasing extension to~$\Xr$.


%

Let us extend $X$ in the same\xx{ way} manner as $\R$ is extended by~(\ref{e:extR}):
$$
\tX=X\cup\xy\{\mainfty,\painfty\},
$$
where $\mainfty$ and $\painfty$ are two\xx{ different} distinct elements\xx{objects} that do not belong to~$X.$
Preorder\xx{Relation} $\succE\!\zaz{\X}\subseteq X\!\times\! X$ is extended to $\tX$ as follows:
$$\succE\!\zaz{\tX}\equiv\;[\succE\!\zaz{\X}\cup\:\{(\painfty,\x)\mid\x\in\tX\}\cup\{(\x,\mainfty)\mid\x\in\tX\}],$$
where $(\painfty,\x)$ and $(\x,\mainfty)$ are pairs of elements of~$\tX.$

Functions $\aa$, $\bb:$\xx{ are defined on $\tX$ from} $\tX\to\tR$ are defined in the same way as in~(\ref{axbx}).

\begin{definition}
\label{def1'}
{\rm
A function $\fP\!:P\to\R,$ where $P\subseteq X,$ is\xx{ said to be} {\em gap-safe increasing with respect to\xx{ the} a preorder $\succEe$} defined on $\X$ (or,\xx{ more} briefly, {\em gap-safe increasing\/})
if $\fPo$ is weakly increasing and 
for any $\x,\,\x'\in\tX,\:$ $\x'\succ \x$ implies $\bb(\x')>\ax$.
}
\end{definition}

The term ``gap-safe increasing'' refers to the property of a function to orderly separate its values (${\bb(\x')>\ax}$) when the corresponding sets of arguments are orderly separated (${\x'\succ\x}$) in $X$; see also Remark~\ref{r:interp}. In~\cite{Cheb02Hagen}, the term ``separably increasing function'' was proposed\xx{, which, however\xx{unfortunately}}, clashing with topological separability, which means the existence of a countable dense subset.

\begin{proposition}
\label{P0}
If $\fPo$ defined on $P\subseteq \X$ is gap-safe increasing\/$,$ then\/$:$\\
$(a)$ $\fPo$ is strictly increasing\/$;$\\
$(b)$ $\fPo$ is\/\footnote{An equivalent formulation is: {\em There is no $\x\in\X$ s.t. $\ax\!=\!+\infty$ or \/$\bx\!=\!-\infty.$
}}
{\xz}upper-bounded on the {\lowercontour} and lower-bounded on the {\uppercontour} of $\x$ for every\/~$\x\in\X.$
\end{proposition}

\xx{Note that there are functions}\xx{It is easy to give an example of a function} It should be noted that there are functions $\fPo$ that are strictly increasing, upper-bounded on all {\lowercontour}s and lower-bounded on all {\uppercontour}s, but are not gap-safe increasing. 
\begin{example}
\label{ex:f}
Consider\xx{ the function}
$$
\fPx=\cases{x\zaz{1},   &$x\zaz{1}\le0,$\cr
            x\zaz{1}-1, &$x\zaz{1}>1,$\cr
}
\label{CoEx} 
$$
where $P=\xx{\,]}(-\infty,0\xy]\,\cup\xx{\,]}(1,+\infty\xx{]})\subset\R^1\xz\equiv\xz\X.$ Function $\fPo$ satisfies (a) and (b) of Proposition~\ref{P0},
but it is not gap-safe increasing. Indeed, $\aa(0)=0=\bb(1).$
\end{example}

\begin{remark}\label{r:interp}
The gap-safe increase of a function can be interpreted as follows.
If $\fPo$ is weakly increasing, then $\x'\succ\x$ implies $\bb(\x')\ge\ax$ for any $\x,\xy\x'\in\X,$ as $(i)$ $\To$ $(iii)$ in Proposition~\ref{P:nonstr}. For the class of strictly increasing functions $\fPo,$ the conclusion cannot be strengthened to $\bb(\x')>\ax,$ as Example~\ref{ex:f} shows.
This stronger conclusion holds\xx{is true1/7} for gap-safe increasing functions, i.e., $\bb(\x')=\ax$ is incompatible for them with $\x'\succ\x.$\xx{ possible only if $\x'\not\succ\x.$} In other words, the \emph{absence of a gap\/} in the values of $\fPo$ between $P$-contours ``{$\x'$ or higher}'' (with infimum given by $\bb(\x')$) and ``{$\x$ or lower}'' (with supremum of $\ax$) \emph{implies\/} $\x'\not\succ\x.$\xx{ (the absence of an) that there is no order gap of the form... (between the corresponding arguments).}
Hence 
the gap-safe increase of a function can be viewed as a kind of continuity of the inverse $f_P^{-1}$ mapping: there is no gap in its values ($\x'\not\succ\x$) whenever\xx{any} there is no gap in the argument\xx{ is lacking} ($\bb(\x')=\ax$).
\end{remark}

\section{A class of extensions of gap-safe increasing functions}
\label{s:main}

Let $\fPo$ defined on any $P\subseteq\X$ be gap-safe increasing.
Theorem~\ref{T} below states that this is a necessary and sufficient condition for the existence of strictly increasing extensions of $\fPo$ to~$\Xr$ provided that $\succE$ enables utility representation. Furthermore, for any\xx{ arbitrary} such a \xx{utility function}representation, the theorem provides an extension of a\xx{ arbitrary} gap-safe increasing function $\fPo$ that combines these two functions.

\smallskip
For any $\a,\b\in\R$ such that $\a<\b,$ let $\uab\!:\X\to\R$ be a utility representation of $\succE$ (i.e., a function strictly increasing w.r.t. $\succEe$\xx{ and bounded}) satisfying
\eq{ge}{ 
\a<\uabx<\b\quad\mbox{for all }\;\x\in\X.
}

For any (unbounded) utility representation of $\succE,$ $u(\x)$, such a function $\uabx$ can be obtained, for example, using transformation
$$
\xx{u\zaz{{\rm example}}(\x)}\uabx=
\frac{\b-\a}{\pi}\left(\arctan\xx{\sum_{i=1}^kx\zaz{i}}u(\x)+\frac{\pi}{2}\right)+\a.
$$

In particular, consider the functions $\u01\!:\X\to\R$ that satisfy
\eq{0<g1<1}{ 
0<\u01(\x)<1.
}

They are normalized versions of the above utilities $\uab$:
\eq{u1}{ 
\u01(\x)=(\b-\a)^{-1}(\uabx-\a),\quad \x\in\X.
}

For any real $\a$ and $\b>\a$ and any utility representations $\u01$ of $\succE,$\xx{that satisfies (\ref{0<g1<1}) and any $\x\in\X$} we define
\begin{eqnarray}
\fx&=&
\max\Big\{\ax,\,\min\big\{\bx,\,\b\big\}-\b+\a\Big\}\big(1-\u01(\x)\big)
\nonumber\\
   &+&\xxy\min\Big\{\bx,\,\max\big\{\ax,\a\big\}-\a+\b\Big\}\,\u01(\x),\quad \x\in\X.
\label{f} 
\end{eqnarray}

\xx{Note that }For an arbitrary gap-safe increasing $\fPo,$ function
$\fo\!:\X\to\R$ given by (\ref{f}) is well defined as the two terms in the right-hand
side\xx{ of (\ref{f})} are finite. This follows from item $(b)$ of Proposition~\ref{P0}.
For preordered sets $\Xr$ that have minimal or maximal elements (see Example~\ref{ex:Nin} in Section~\ref{s:Related}, where $\Xr$ has a maximal element), this is ensured by introducing\xx{using} the augmented sets $\tX$ in the definition\xx{~\ref{def1'}} of a gap-safe increasing function. Indeed, since\xx{ for every $\fP$,} $\bb(\painfty)=+\infty,$ ${\aa(\mainfty)=-\infty,}$ and $\painfty\succ\x\succ\mainfty$ for all $\x\in\X,$ Definition~\ref{def1'} provides $\bb(\painfty)>\ax$ and $\bx>\aa(\mainfty),$ hence $+\infty>\ax$ and ${\bx>-\infty,}$\xx{ implies that for a gap-safe increasing function,}\xx{ whenever $x\in X$ is a maximal and a minimal element of $\succE$ on $X,$ respectively.} i.e., $\fP$ is upper-bounded on all {\lowercontour}s and lower-bounded on all {\uppercontour}s, ensuring the correctness of definition~(\ref{f}). If $\Xr$ has neither minimal nor maximal elements (like the Pareto preorder\xx{ $\succE$} on~$\Rk$), then the replacement of $\tX$ with $X$ in Definition~\ref{def1'} does not alter the class of gap-safe increasing functions. 

We now formulate the main result\xx{ of this paper}.

\begin{theorem}
\label{T}
Suppose that a preorder $\succE$ defined on $X$ has a utility representation and
$\fPo$ is a real-valued function defined on some ${P\subseteq\X}$.
Then $\fP$ is strictly monotonically extendable to $\Xr$ if and only if $\fP$ is gap-safe increasing.

Under these conditions\/$,$
function $\fo$ defined by $(\ref{f}),$ where $\u01$ is any utility representation of\/ $\succE$ that satisfies~$(\ref{0<g1<1})$ and $\a<\b,$ is a strictly increasing extension of $\fPo$ to\/~$\Xr$.
\end{theorem}

\section{Extension of utility: Additional representations}
\label{s:alter}

The class of extensions introduced by Theorem~\ref{T} allows alternative representations that clarify its properties.
They are given by\xx{ the following} Propositions~\ref{P1}--\ref{P3}.
\begin{proposition}
\label{P1}
If\/ $\uab\!:\X\to\R$ is a utility representation of\/ $\succE$ satisfying~$(\ref{ge})$ and $\fPo\!:$ $P\to\R,$ where $P\subseteq\X,$ is gap-safe increasing\/$,$ then
\begin{eqnarray}
\!\!\!\!\fx=(\b-\a)^{-1}\!\!
&\Big(&\!\max\Big\{\ax-\a,\,\xxy\min\big\{\bx-\b,\,0\big\}\!\Big\}\big(\b-\uabx\big)
\nonumber\\
&+&\min\Big\{\bx-\b,\,\max\big\{\ax-\a,\,0\big\}\!\Big\}\big(\uabx-\a\big)\!\Big)
\nonumber\\
&+&\uabx
\label{f'}
\end{eqnarray}
is a strictly increasing extension of $\fPo$ to\/~$\Xr,$ and $\fx$ coincides with function\xx{ introduced by}~$(\ref{f}),$ where $\u01$ is\xx{ defined} related to $\uab$ by\/~$(\ref{u1})$.
\end{proposition}

The order of proofs in Section~\ref{a:Proofs} is as follows.
Verification of the second statement of Proposition~\ref{P1} is straightforward and\xx{ it} is omitted.
\xx{It is straightforward to verify that function (\ref{f'}) coincides with (\ref{f}) where $\u01$ is defined by~(\ref{u1}).}
This statement is used\xx{ in Section~\ref{a:Proofs}} to prove Proposition~\ref{P3}, which implies Proposition~\ref{P2}, and they both are\xx{ applied} used\xx{employed} in the proof of Theorem~\ref{T}, which in turn implies the first statement of Proposition~\ref{P1}.
%
%

\smallskip
To simplify (\ref{f'}), we
partition $\X\setminus P$ into four regions determined by $\succE$ and~$P$:
\eq{ALUN}{ 
\begin{array}{rrr}
A\,=\,\big\{\x\in \X\setminus P\;\,\big|\msd&
         \Pmx\ne\emptyset\;\;\mbox{and\,}&\Ppx\ne\emptyset\big\},\\
L\,=\,\big\{\x\in \X\setminus P\;\,\big|\msd&
     \Pmx=\emptyset\;\;\mbox{and\,}&\Ppx\ne\emptyset\big\},\\
U\,=\,\big\{\x\in \X\setminus P\;\,\big|\msd&
         \Pmx\ne\emptyset\;\;\mbox{and\,}&\Ppx=\emptyset\big\},\\
N\,=\,\big\{\x\in \X\setminus P\;\,\big|\msd&
     \Pmx=\emptyset\;\;\mbox{and\,}&\Ppx=\emptyset\big\}.
\end{array}
}

\xx{ Obviously}Clearly these regions are pairwise disjoint and\/ $\X=P\cup A\cup L\cup U\cup N$.

\begin{proposition}
\label{P2}
If\/ $\uab\!:\X\to\R$ is a utility representation of\/ $\succE$ satisfying\/~$(\ref{ge})$ and $\fPo\!:$ $P\to\R,$ where $P\subseteq\X,$ is gap-safe increasing\/$,$ then
function $\fo$ defined by $(\ref{f'})$ can be represented as follows\/$:$
\eq{f''}{ 
\fx=\cases{\fPx,                              &$\x\in P,$\cr
       \xxy\min\big\{\bx -\b,\,0\big\}+\uabx, &$\x\in L,$\cr
           \max\big\{\axu-\a,\,0\big\}+\uabx, &$\x\in U,$\cr
           \uabx,                             &$\x\in N,$\cr
           \mbox{\xx{not simplified }expression $(\ref{f'}),$}
                                              &$\x\in A$.\cr
}}
\end{proposition}

Proposition~\ref{P2}\xx{ clarifies} highlights the role of $\uab$ in~(\ref{f'}).\xx{of the extension~$\fo$.} \xx{According to (\ref{f''}),}Function $\fo$ reduces to $\fP$ on $P$ and to $\uab$ on $N$ whose elements are {$\succE$-incomparable} with those of~$P$. Moreover, $\fx=\uabx$ on the part of $L$ where $\bx\ge\b$ and on the part of $U$ where $\ax\le\a.$ On the complement parts of $L$ and $U$, $\fx=\bx+(\uabx-\b)$ and $\fx=\ax+(\uabx-\a),$ respectively. On $A,$ (\ref{f'}) is not simplified. This fact and the ambiguity on\xx{in} $L$ and $U$ prompt us to make another\xx{ a different} partition of~$X.$

Consider four regions that depend on $\succE,$ $P,$ $\fP,$ $\a,$ and $\b$:
\eq{e:Sdom}{
\begin{array}{lll}
S_1&=&\big\{\x\in \X\:\big|\;\bx-\ax\le\b-\a\big\},\cr
S_2&=&\big\{\x\in \X\:\big|\;\bx-\ax\ge\b-\a\mbox{\, and }\,\bx\le\b\big\},\cr
S_3&=&\big\{\x\in \X\:\big|\;\bx-\ax\ge\b-\a\mbox{\, and }\,\ax\ge\a\big\},\cr
S_4&=&\big\{\x\in \X\:\big|\;\ax\le\a\mbox{\, and }\,\bx\ge\b\big\}.
\end{array}
}
It is easily seen that $\X=S_1\cup S_2\cup S_3\cup S_4,$ whereas the $S_i$-regions are not disjoint.
This\xx{ partition} decomposition allows us to express\xx{yields an expression for} $\fx$ without $\min$ and~$\max$.

\begin{proposition}
\label{P3}
For a gap-safe increasing $\fPo,$ $\fo$ defined by $(\ref{f})\!$\xx{ or $(\ref{f'})$} can be represented as follows\/$,$ where $\u01$ and $\uab$ are\xx{ utility} representations of\/ $\succE$ related by~$(\ref{u1})\!:$
\eq{f'''}{ 
\fx=\cases{\ax\big(1-\u01(\x)\big)+\bx\,\u01(\x),   &$\x\in S_1,$\cr
             \bx+\uabx-\b,                          &$\x\in S_2,$\cr
            \axu+\uabx-\a,                          &$\x\in S_3,$\cr
             \uabx,                                 &$\x\in S_4.$\cr
}
}
\end{proposition}

Thus, on $S_1,$ $\fx$ is a convex combination of $\bxl$ and $\axl$ with coefficients $\u01(\x)$ and $(1-\u01(\x)),$ respectively.
The regions $S_1,S_2,S_3,$ and $S_4$\xx{ meet} intersect on some parts of the border sets ${\bx-\ax}=\b-\a$, $\ax=\a$, and $\bx=\b$. Accordingly, the expressions of $\fo$ given by\xx{in} Proposition~\ref{P3} are concordant on these intersections.

\begin{corollary}
\label{c:reduce}
In the notation and assumptions of Proposition~$\ref{P3},$ $N\subseteq S_4.$\\
For any $\x\in X,$ $\axl=\bxl$ implies $\fx=\ax.$ In particular$,$
if\/ $\x\approx\p$ for some $\p\in P,$ then $\fx=\fPp$ and\/ $\x\in S_1.$
\end{corollary}

\section{\xx{Extendability}Extension of functions defined on Pareto sets}
\label{SPar}

Consider the case where $P$ is a Pareto set.
Such a set comprises elements that are mutually undominated.

\begin{definition}
{\em
A subset $P\subseteq\X$ is called a {\em Pareto set\/} in $\Xr$ if there are no $\p,\xy\p'\in P$ such that $\p'\succ\p,$ where $\succ$ is the asymmetric part of~$\succE$.
}
\end{definition}

For functions defined on Pareto sets $P$, the necessary and sufficient condition of extendability to $\Xr$ given by Theorem~\ref{T} reduces to the\xx{ condition of} boundedness on all $P$-contours (which appeared in\xx{ item $(b)$ of} Proposition~\ref{P0}) supplemented by condition (\ref{e:appro}):
$\eqeqA$.

\begin{lemma}\label{l:sin-contour}
A function $\fPo$ defined on a Pareto set $P\subseteq\X$ is gap-safe increasing with respect to a preorder $\succEe$ defined on $\X$ if and only if $\fPo$ is upper-bounded on all {\lowercontour}s\/$,$ lower-bounded on all {\uppercontour}s\/$,$ and satisfies\xx{ condition $(\ref{e:appro})${\rm:}}
$\eqeqA,$
where $\approx$ is the symmetric part of\/~$\succE$.
\end{lemma}

By the transitivity of $\succE$\xx{ $\succN$ and $\approx$}, for any Pareto set $P,$ the sets $P\cup A$ and $S_1$ have a simple structure described in the following lemma.\xx{ the form $A\,=\,\big\{\x\in \X\setminus P\;\big|\;\x\approx\p\/\;\mbox{for some}\:\p\in P\big\}$.}

\begin{lemma}
\label{l:setA}
Under the conditions of Lemma~$\ref{l:sin-contour},$
$S_1=P\cup A=\{\x\in\X\mid\exists\xy\p\in P\!:\p\approx\x\},$
where $S_1$ and $A$ are defined by $(\ref{e:Sdom})\xz$ and $(\ref{ALUN}),$ respectively.
\end{lemma}

Lemmas~\ref{l:sin-contour} and~\ref{l:setA}, Propositions~\ref{P2} and~\ref{P3}, and Corollary~\ref{c:reduce} provide the following special case of Theorem~\ref{T} for Pareto sets.

\begin{corollary}
\label{c:Par}
Suppose that a preorder $\succE$ on $X$ has a utility representation $\uab$\xx{$\!:\X\to\R$} satisfying $(\ref{ge})$ and $P\subseteq\X$ is a Pareto set.
Then a function $\fPo\!:P\to\R$ is strictly monotonically extendable to $\Xr$ if and only if\xx{ $\fPo$} it is upper-bounded on all {\lowercontour}s\/$,$ lower-bounded on all {\uppercontour}s\/$,$ and satisfies 
$\eqeqA,$
where $\approx$ is the symmetric part of\/~$\succE$.

Under these conditions\/$,$ the function $f\!:\X\to\R$ such that\/$:$

\begin{tabular}{ll}
$\quad\fx=\fPp,$                                          &\;\;whenever $\p\approx\x$ and\xx{for some} $\p\in P;$\\
$\quad\fx$ is defined by $(\ref{f''})\xz$ or $(\ref{f'''}),$ &\;\;\xx{otherwise\/$,$ i.e.$,$}when $\x\not\in P\cup A=S_1$
\end{tabular}

\noindent
is a strictly increasing extension of\/ $\fPo$ to $\Xr$ coinciding with\/~$(\ref{f'}).$
\end{corollary}

It follows from Corollary~\ref{c:Par} that for a Pareto set $P$, functions $\fP$ and $\uab$ influence $f$ almost symmetrically: $f$ reduces to $\fP$ on $P\cup A=S_1,$ to $\uab$ on $S_4$\xx{$\supseteq N$}, and is determined by the sum $\bx+\uabx$ or $\axu+\uabx$ on $S_2\cup S_3.$

\smallskip
Results\xx{ closely} related to Theorem~\ref{T} and Corollary~\ref{c:Par} were used in \cite{CheSha98AOR,Che23JCNCentr} to construct implicit forms of\xx{ monotonic monotone} scoring procedures for preference aggregation and evaluation of the centrality of nodes.
More specifically, theorems\xx{results} of this type\xx{kind make it possible} allow us to move from axioms that determine a positive impact of the comparative results of objects on their functional scores to the conclusion that the scores satisfy a system of equations determined by a strictly increasing function.

\section{Connections to related work}
\label{s:Related}

Problems of extending real-valued functions\xx{ from a subset to a set} while preserving monotonnicity\xx{certain properties} (\xx{also}sometimes called lifting problems) have been considered primarily in topology. Therefore, continuity was usually\xx{mainly considered as} a\xx{ necessary} property to be preserved. This strand of literature started with the following theorem of general topology.

\TheoCite{Urysohn's extension theorem~\cite{Urysohn25Extension}}
{A topological space $(X,\tau)$ is normal\xy\footnote{A topological space $(X,\tau)$ is called {\em normal\/} if for any two disjoint closed subsets of $X$ there are two disjoint open subsets each covering one of the closed subsets.}
if and only if every continuous real-valued function $\fP$ whose domain is a closed subset ${P\subset X}$  can be extended to a\xx{ real-valued} function\xx{ $f$} continuous on~$X.$}

For metric spaces, a counterpart of this theorem was proved by Tietze~\cite{Tietze1915Funktionen}.

Nachbin \cite{Nachbin65TopolOrderNo}\xx{ proved} obtained extension theorems for functions defined on preordered spaces. In his terminology,\xx{ He introduces the notion of a normally ordered space. He says that} a topological space $\XX$ equipped with a preorder $\succE$ is {\em normally preordered\/} if for any two disjoint closed sets $F_0,\,F_1\subset X,$ $F_0$ being {\em decreasing\/} (i.e., with every $\x\in F_0$ containing all $\y\in X$ such that $\y\precE\x$) and $F_1$ {\em increasing\/} (with every $\x\in F_1$ containing all $\y\in X$ such that $\y\succE\x$), there exist disjoint open sets $V_0$ and $V_1,$ decreasing and increasing respectively, such that $F_0\subseteq V_0$ and $F_1\subseteq V_1.$ The space is {\em normally ordered\/} if, in addition, its preorder $\succE$ is antisymmetric (i.e., is a partial order).

\TheoCite{Nachbin's lifting theorem \cite{Nachbin65TopolOrderNo}\xx{[Theorem~6 in Chapter~1]} for compact sets in ordered spaces}{In any normally ordered space $\XX$ whose partial order $\succE$ is a closed subset of\/ ${X\!\times\! X,}$ every continuous weakly increasing real-valued function defined on any compact set ${P\subset X}$ can be extended to $X$ in such a way as to remain continuous and weakly increasing.}

An analogous theorem for more general\xx{ case of a} normally {\em preordered\/} spaces is~\cite[Theorem~3.4]{Minguzzi13NormPre}.
Sufficient conditions for\xx{ a topological space} $\XX$ to be normally preordered are: (a)~compactness of $X$ and belonging of $\succE$ to the class of\xx{ equipped with a} closed\xx{??} partial orders \cite[Theorem~4 in Chapter~1]{Nachbin65TopolOrderNo} (this result was strengthened in~\cite{Minguzzi13NormPre}); (b)~\xx{completeness}connectedness and closedness of $\succE$~\cite{Mehta77TopoOrd}.

\xx{Some other}Additional utility extension theorems in which $P$ is a compact set, $\fP$ is continuous, and $f$ is required to be continuous and weakly\xx{ (or strictly)} increasing as well as~$\fP$ are discussed in~\cite{Evren21Extension}.

The extendability of continuous functions defined on\xx{ generic} non-compact sets $P$\xx{ of topological spaces} requires a stronger condition.
It can be formulated as follows.

For a function $\fP\!:P\to\R,$ where $P\subseteq X,$ let
the {\em lower $\fP$-contour\/} and
the {\em upper $\fP$-contour of\/ ${r\in\R}$}\xx{ be} denote the sets
$\Pmr   \equiv\{\p\in P\mid\fP(\p)\le r\} $ and
$\Ppr{r}\equiv\{\p\in P\mid\fP(\p)\ge r\},$
respectively.
Let us\xx{We} say that $\fP$ is {\em inversely closure-increasing\/} if for any $r,\,r'\in\R$ such that ${r<r',}$ there exist two disjoint closed subsets of $X$: a decreasing set containing $\Pmr$ and an increasing set containing $\Ppr{r'}.$
%

\TheoCite{Nachbin's lifting theorem \cite{Nachbin65TopolOrderNo}\xx{[Theorem~2 in Chapter~1]} for closed sets in preordered spaces}
{In any normally preordered space $\XX,$ a continuous weakly increasing bounded function $\fP$ defined on a closed subset $P\subset X$ can be extended to $X$ in such a way as to remain continuous\/$,$ weakly increasing\/$,$ and bounded if and only if $\fP$ is inversely closure-increasing.}

For several other results regarding the extension of weakly increasing functions defined on non-compact sets $P,$ we refer to \cite{Herden89Liftings,Herden90LiftingN,Minguzzi13NormPre}.

Theorems on the extension of {\em strictly\/} increasing functions were obtained in~\cite{Herden89Liftings,Husseinov10Monotonic,Husseinov18ExtStrict,Husseinov21Extension}. Herden's\xx{ general} Theorem~3.2 \cite{Herden89Liftings} contains a compound condition consisting of several arithmetic and set-theoretic parts, which is not easy to grasp\xx{ verify or interpret}. To formulate a more transparent result \cite[Theorem~2.1]{Husseinov18ExtStrict}, let us introduce the following notation. Using Definition~\ref{def23}, for any $Z\subseteq X$ define the {\em decreasing cover of~$Z,$} $d(Z)=\bigcup_{\z\in Z\msd}\Xm(\z)$ and the {\em increasing cover of~$Z,$} $i(Z)=\bigcup_{\z\in Z}\Xp(\z)$.
In these terms, $Z$ is decreasing (increasing) whenever $Z=d(Z)$ (resp., ${Z=i(Z)}$). 
A preorder\xx{ $\succE$} is said to be {\em continuous\/} \cite{Mccartan71ContinPreord} if for every open $V\subset\X,$ both $d(V)$ and $i(V)$ are open.
A~preorder $\succE$ is {\em separable\/}\footnote{On connections between versions of preorders' separability and denseness, see~\cite{Herden89Existence}.} if there exists a countable ${Z\subseteq X}$ such that $[\x,\,\x'\in\X\mbox{ and }\x\prec\x']\To[{\x,\,\x'\in Z}\mbox{ or }(\x\prec\z\prec\x'\mbox{ for some }\z\in Z)].$
For ${\x\in X}$ denote by $\Vx{d}$ and $\Vx{i}$ the collections of open decreasing and open increasing sets containing $\x$, respectively.

\TheoCite{H{\"u}sseinov's extension theorem \cite{Husseinov18ExtStrict} for strictly increasing functions}{In any normally preordered space $\XX$ with a separable and continuous preorder~$\succE,$ a continuous strictly increasing function $\fP$ defined on a nonempty closed subset $P\subset X$ can be extended to $X$ in such a way as to remain continuous and strictly increasing if and only if $\fP$ is such that for any $\x,\,\x'\in\X,\:$ $\x'\succ\x$ implies $\bb(\x')>\ax$
and for any $\x\in\X,\xy$ $M(\x)\ge m(\x),$ where

\smallskip
\centerline{$\displaystyle
m(\x)=\!\inf_{V_d\in\Vx{d}}\!\!\sup\{f(\p)\!\mid\xz\p\xz\in\xz P\cap V_d\}\mbox{ and\/ }
M(\x)=\!\sup_{V_i\in\Vx{i}}\!\!\inf\{f(\p)\!\mid\xz\p\xz\in\xz P\cap V_i\}
$}

\smallskip\noindent
with the convention\xx{agreement} that
$m(\x)=\inf\{f(\p)\!\mid\xz\p\xz\in\xz P\}$
if $P\cap V_d=\emptyset$ for some $V_d\in\Vx{d}$ and\/
$M(\x)=\sup\{f(\p)\!\mid\xz\p\xz\in\xz P\}$
if $P\cap V_i=\emptyset$ for some $V_i\in\Vx{i}.$
}

This theorem is a topological counterpart of the first part of our Theorem~\ref{T}. 
Consider the discrete topology in which every subset of $\X$ is open. Then the space $\XX$ is normally preordered and the preorder $\succE$ is continuous, as well as any function~$\fP$.
The separability of $\succE$ in H{\"u}sseinov's theorem ensures its representability by utility, which is explicitly assumed in Theorem~\ref{T}.

Condition $M(\x)\ge m(\x)$\xx{ for all $\x\in X$} reduces to $\bdx\ge\adx\xx{\,\mbox{for all }\;\x\in\X},$ 
where $\bdx$ and $\adx$ modify $\bx$ and $\ax$ by taking values $\sup\{\fP(\p)\!\mid\xz\p\xz\in\xz P\}$ or $\inf\{\fP(\p)\!\mid\xz\p\xz\in\xz P\}$ instead of $+\infty$ or $-\infty$ when $\Ppx=\emptyset$ or $\Pmx=\emptyset,$ respectively. It is easily seen that conditions $\bdx\ge\adx$ and $\bx\ge\ax$ are equivalent (cf.\ Remark~\ref{r:sup}), therefore, by $(i)\!\ToTo\!(ii)$ of Proposition~\ref{P:nonstr}, $M(\x)\ge m(\x)$ for all $\x\in X$ reduces in the discrete topology to the weak increase of~$\fP$.

The last condition, $\x'\succ\x\To\bb(\x')>\ax,$ proposed in \cite{Cheb02Hagen}, is required for all $\x,\x'\in\X$ in the above theorem and for all $\x,\x'\in\tX$ in Theorem~\ref{T} (forming, by Definition~\ref{def1'}, the main part of gap-safe increase). 
\xx{It is the only}This difference is significant. Let us illustrate it with the following example. 
\begin{example}
\label{ex:Nin}
$\X=\Z\setminus\N=\{\bm{0},\bm{-1},\bm{-2},...\};$ $\succE\:=\!\bigcup_{\x\in X\setminus\{\bm0\}}\{(\bm{0},\x),(\x,\x)\}\cup\{(\bm0,\bm0)\};$ $P=\X\setminus\{\bm{0}\};$ $\fPp=-\p$ for all $\p\in P$.

Then $\fP$ has no strictly increasing extension to\/ $\Xr$ and is not gap-safe increasing, since $\painfty\succ\bm{0}$, but $+\infty=\bb(\painfty)\not>\aa(\bm{0})=+\infty.$ However, $\x'\succ\x\To\bb(\x')>\ax$ for all $\x,\x'\in\X,$ therefore,\xx{ and, as all conditions of the above theorem are satisfied, it} the above theorem claims that $\fP$ is strictly monotonically extendable to~$\Xr$.
\end{example}

The reason is that \cite[Theorem~2.1]{Husseinov18ExtStrict} was actually proved for a bounded function $\fP$, however, the boundedness condition was removed by a remark erroneously claiming that this condition was\xx{is} not essential.
The method of extension proposed in the present paper differs from the classical approach, which is systematically applied to continuous functions.

In \cite{Husseinov10Monotonic}, H\"{u}sseinov shows that condition $M(\x)\ge m(\x)$ for all $\x\in X$ is equivalent to the necessary and sufficient extendability condition for a weakly increasing bounded function $\fP$ defined on a closed subset of a preordered space, i.e., to the Nachbin property of being inversely closure-increasing.

The problem of extending utility functions without continuity constraints was considered in~\cite{Cheb02Hagen} with the focus on \xx{. The main content of this paper relates to}the functions representing Pareto {\em partial\/} orders on Euclidean spaces. Partial orders are antisymmetric preorders,\xx{ as a more general class,} therefore, preorders are more flexible allowing symmetry ($\x\succE\y,\,\y\succE\x$) on a pair of distinct elements, while partial orders only allow ``negative'' ($\x\not\succE\y,\,\y\not\succE\x$) symmetry. Symmetry is an adequate model of equivalence between objects (which suggests the same value of the utility function), while ``negative'' symmetry can model\xx{corresponds to (matches)} the absence of information, which is generally compatible with unequal utility values.

\section{Conclusion}
\label{s:Conclu}

The paper presents a strict-extendability condition and\xx{ an algorithm for extending} a wide class of extensions of utility functions $\fP$ defined on an arbitrary subset $P$ of an arbitrary set $X$ equipped with a preorder.
It can be observed that the key condition of gap-safe increase of $\fP$ has a similar structure as that of inverse closure-increase, which is equivalent to the extendability of a continuous weakly increasing function $\fP$ defined on a closed subset $P\subset\X$ (see~\cite{BridgesMehta95ReprPrefOrd} for a related discussion). Moreover, as mentioned in Section~\ref{s:Related}, the latter ``inverse'' condition has an equivalent ``direct'' counterpart. Relationships of this kind deserve further study.

Among other problems, we\xx{ would like to} mention: (1)~finding relationships between various extensions \xx{developed}proposed earlier for continuous functions and the class of extensions (\ref{f}) described in Theorem~\ref{T}; (2)~characterizing the complete class of extensions of $\fP$ to $\Xr$ (and, for example, to~$(\Rk,\mbox{Pareto preorder}$)); (3)~exploring the extension problem with $\R$ as the range of\xx{ $\fP$ and} $f$ replaced by certain other posets. \xx{(partially ordered spaces).}

\section{Proofs}
\label{a:Proofs}


{\em Proof of Proposition~\ref{P:nonstr}.}
$(i)\!\To\!(ii)$. Let $(i)$ hold. 
For any $\x\!\in\!\xx{\tX}X,$ if ${\Pmx=\emptyset}$ or $\Ppx=\emptyset,$ then $\ax=-\infty$ or $\bx=+\infty,$ respectively, with $\bx\ge\ax$ in both cases. Otherwise, $\p'\in\Pmx$ and $\p''\in\Ppx$ imply $\p''\succE\x\succE\p',$ and $\p''\succE\p'$ by the transitivity of~$\succEe$. Hence $\fP(\p'')\ge\fP(\p')$ by~$(i).$\xx{ as $\fP$ is weakly increasing.} Therefore, $\inf\{\fP(\p'')\mid\p''\in\Ppx\}\ge\sup\{\fP(\p')\mid\p'\in\Pmx\},$ i.e., $\bx\ge\ax.$
%

$(ii)\!\To\!(iii)$. Let $(ii)$ hold. Then for any $\x,\x'\in\xx{\tX}X$ such that $\x'\succEe\x$ using (\ref{e:abmono}) we get $\bb(\x')\ge\bx\ge\ax.$

$(iii)\!\To\!(ii)$. As $\succE$ is reflexive, $(ii)$ is a special case of~$(iii)$.

%
$(iv)\ToTo(i)\ToTo(v)$.
[For all $\,\p\in P,$ $\fP(\p)\ge \ap$] $\ToTo$
[for all $\,\p,\,\p'\in P,$ $(\p\succEe\p')\To(\fP(\p)\ge\fP(\p'))$] $\ToTo$
[for all $\,\p'\in P,$ $\bb(\p')\ge\fP(\p')$].

$(vi)\To(iv)$. [$(vi)$ and the last inequality of (\ref{a>f>b})] $\To$ $(iv).$

$(ii)\To(vi)$ as $(vi)$ is a special case of $(ii)$.
\qed

\medskip
{\em Proof of Proposition~\ref{P0}.}
Let $\fPo$\xx{ $P\subseteq X,$} be gap-safe increasing.

$(a)$ Assume that $\fPo$ is not strictly increasing. Since $\fP$\xx{ satisfies $(iii)$ of Proposition~\ref{P:nonstr}, it} is weakly increasing,\xx{ this implies that} there are $\p,\,\p'\in P$ such that $\p'\succ\p$ and $\fPp=\fP(\p')$. Then, by (\ref{a>f>b}), $\ap\ge\fPp=\fP(\p')\ge \bb(\p')$ holds$,$ i.e., $\fPo$ is not gap-safe increasing. Therefore, the assumption is wrong.

$(b)$ Let $\Pmx$ be the \lowercontour\ of some $\x\in X$.\xx{ Consider any $\x'\in\tX$ such that $\x'\succ\x$.} By definition, $\painfty\in\tX$ and $\painfty\succ\x.$ Since $\fPo$ is gap-safe increasing,\xx{ $\bb(\x')>\ax$} $+\infty=\bb(\painfty)>\ax$.\xx{ Therefore, $\ax<+\infty$.} Since $\ax=\sup\big\{\fP(\p)\mid\p\in \Pmx\big\}$, $\fPo$ is upper-bounded on~$\Pmx$. Similarly, $\fPo$ is lower-bounded on all {\uppercontour}s.
%
%
\xx{This completes the proof.}
\qed

\medskip
Next we prove Proposition~\ref{P3}; then it will be used to prove Proposition~\ref{P2} and Theorem~\ref{T}.

\medskip
{\em Proof of Proposition~\ref{P3}.}
Let $\x\in S_1.$ Since $\bx-\ax\le\b-\a,$ we have
\begin{eqnarray*}
\min\big\{\bx,\,\b\big\}-\ax &\,\le\,& \b-\a,\\
\bx-\max\big\{\ax,\a\big\}   &\,\le\,& \b-\a,
\end{eqnarray*}\vspace{-1em}
hence

\vspace{-1em}
\begin{eqnarray*}
\ax&\,\ge\,&\min\big\{\bx,\,\b\big\}-\b+\a,\\
\bx&\,\le\,&\max\big\{\ax,\a\big\}-\a+\b.
\end{eqnarray*}
Therefore, (\ref{f'}) reduces to $\fx=\ax\big(1-\u01(\x)\big)+\bx \u01(\x)$.

Let $\x\in S_2.$ Inequalities $\bx-\ax\ge\b-\a$ and $\bx\le\b$ imply
$\ax\le\a$, hence (\ref{f}) reduces to $\fx=\bx+\uabx-\b$.

Let $\x\in S_3.$ Inequalities $\bx-\ax\ge\b-\a$ and $\ax\ge\a$ imply
$\bx\ge\b$, hence (\ref{f}) reduces to $\fx=\ax+\uabx-\a$.

\xx{The proof for the case of $\x\in S_4$ is straightforward.}
Finally, let $\x\in S_4,$ i.e., $\ax\le\a$ and $\bx\ge\b.$
Substituting $\max\xz\big\{\ax-\a,\,0\big\}=0$ and $\min\xz\big\{\bx-\b,\,0\big\}=0$ into (\ref{f'}) yields $\fx=\uabx.$
\qed

\medskip
{\em Proof of Proposition~\ref{P2}.}
Let $\x\in P.$ Then\xx{, by item (d) of Proposition~\ref{P0}, $\bx=\ax=\fPx$, therefore,} by (\ref{a>f>b}) and (\ref{ge}), $\bx-\ax\le\b-\a$, hence  
$\x\in S_1$. Using Proposition~\ref{P3}, we have $\fx=\fPx\big(1-\u01(\x)\big)+\fPx\, \u01(\x)=\fPx$.

Let $\x\in U$. Then $\bx=+\infty$, hence (\ref{f'}) reduces to $\fx=\max\big\{\ax-\a,\,0\big\}+\uabx$. Similarly, if  $\x\in L$, then $\ax=-\infty,$ and (\ref{f'}) reduces to $\fx=\min\big\{\bx-\b,\,0\big\}+\uabx$.

Finally, if $\x\in N,$ then $\ax=-\infty$ and $\bx=+\infty$, whence $\ax<\a$ and $\bx>\b$, and Proposition~\ref{P3} provides $\fx=\uabx$.
\qed

\medskip
{\em Proof of Theorem~\ref{T}.}
Suppose that $\fPo$ is strictly monotonically extendable to $\Xr$. Then $\fPo$ is strictly increasing w.r.t.~$\succE$.\xx{ and therefore weakly increasing. Then, by Proposition~\ref{P:nonstr}, $\x,\x'\in \X$ and $\x'\succEe\x$ imply that $ \bb(\x')\ge\ax.$}
Assume that $\fPo$ is not gap-safe increasing. This implies that there are $\x,\,\x'\in\tX$ such that $\x'\succ\x$ and $\bb(\x')\le\ax.$
If $\x,\,\x'\in\X,$ then using this inequality, the definition of $\aa$ and $\bb,$ and the strict monotonicity of $\fo,$\xx{ fact that $\fo$ is strictly increasing} we obtain $f(\x')\le \bb(\x')\le\ax\le\fx,$ whence $f(\x')\le\fx$, and as $\x'\succ\x,$ $\fo$ is not strictly increasing.
%
Therefore,\xx{ $\x\in\tX\setminus X$ or $\x'\in\tX\setminus X.$} $\{\x,\,\x'\}\not\subseteq X.$
If $\x\in\tX\setminus X,$ then $\x'\succ\x$ implies $\x=\mainfty$ and $\x'\in X\cup\{\painfty\}.$ By the assumption, $\bb(\x')\le \ax=\sup\emptyset=-\infty,$ hence $\bb(\x')=-\infty$, thus, $\x'\ne\painfty$ and $\x'\in X.$ Since, $\bb(\x')=-\infty$, $f(\x')$ cannot be assigned a value compatible with the strict monotonicity of $\fo$, whence $\fP$ is not strictly monotonically extendable to~$\Xr$, a contradiction. The case of $\x'\in\tX\setminus X$ is considered similarly. It is proved that $\fPo$ is gap-safe increasing whenever $\fPo$ is strictly monotonically extendable to~$\Xr$.

Now let $\fPo$ be gap-safe increasing. By Proposition~\ref{P2}, the restriction of $\fo$ to $P$ coincides with~$\fPo$.

\xx{Let us}It remains to prove that $\fo$ is strictly increasing on~$\X$. This can be\xx{ demonstrated} shown directly by analyzing\xx{ considering(decomposing)} expression~(\ref{f}). Here, we give\xx{ another} a\xx{ different} proof that does not\xx{ require additional calculations} require the analysis of special cases with $\xy\min\xy$ and $\xy\max$.

By Proposition~\ref{P3}, function~(\ref{f}) coincides with (\ref{f'''}), where $\uab$ and $\u01$ are related by~(\ref{u1}).

We will use Lemma~\ref{l:strict}.
First, consider any $\x,\,\x'\in\X$ such that $\x'\approx\x$ and show that $f(\x')=\fx.$
By (\ref{e:ApprCont}), $\aa(\x')=\ax$ and  $\bb(\x')=\bx.$ Furthermore, $\uab$ and $\u01$ are strictly increasing with respect to $\succE$ by definition, hence $\uab(\x')=\uabx$ and $\u01(\x')=\u01(\x).$ Therefore, by (\ref{f'''}), $f(\x')=\fx$ holds.

Now suppose that $\x,\,\x'\in\X$ and $\x'\succ\x$. Then, by (\ref{e:abmono}) and the strict monotonicity of $\uab$ and $\u01$, we have
\eq{pre}{ 
\begin{array}{rcl}
\uab(\x')&\,  >\,&\uabx,   \\
\u01(\x')&\,  >\,&\u01(\x),\\
\aa(\x') &\,\ge\,&\ax,     \\
\bb(\x') &\,\ge\,&\bx.
\end{array}
}

\xx{Suppose{ first} that}Let $\x$ and $\x'$ belong to the same region:
$S_2,S_3,$ or $S_4$. \xx{Then}Ineqs (\ref{pre}) yield
\begin{eqnarray}\nonumber
\bb(\x')+\uab(\x')-\b&\,>\,&\bx+\uabx-\b, \cr
\aa(\x')+\uab(\x')-\a&\,>\,&\ax+\uabx-\a, 
\label{mo} 
\end{eqnarray}
hence, by (\ref{f'''}), $\fo$ is strictly increasing on each of these regions.

If $\x,\x'\in S_1$, then by (\ref{f'''}), (\ref{pre}), (\ref{0<g1<1}), and item $(ii)$~of Proposition~\ref{P:nonstr},
\begin{eqnarray}
f(\x')-\fx
&\ge& \ax\big(1-\u01(\x')\big)+\bx\xy\u01(\x')\nonumber\\
&-&   \ax\big(1-\u01(\x )\big)-\bx\xy\u01(\x )\nonumber\\\nonumber
&=&   \big(\bx-\ax\big)\big(\u01(\x')-\u01(\x)\big)\ge0.
\end{eqnarray}

This implies that $f(\x')=\fx$ is possible only if $\bb(\x')=\bx$ and $\bx=\ax$, hence only if $\bb(\x')=\ax$. The last equality is impossible, since $\fPo$ is gap-safe increasing by assumption. Therefore, $f(\x')>\fx$, and $\fo$ is strictly increasing on~$S_1$.

Let now $\x$ and $\x'$ belong to different regions $S_i$ and~$S_j$.
Consider the points that represent $\x$ and $\x'$ in the 3-dimensional space with axes corresponding to $\aa(\cdot)$, $\bb(\cdot)$, and $\u01(\cdot)$. Let us connect these points, $\big(\ax,\bx,\u01(\x)\big)$ and $\big(\aa(\x'),\bb(\x'),\u01(\x')\big)$, by a line segment. The projections of this segment and the borders of the regions $S_1,S_2,S_3,$ and $S_4$ onto the plane $\u01=0$ are illustrated in Fig.~\ref{F1}.

\begin{figure}[ht]
\unitlength 0.90mm
\linethickness{0.3pt}
\begin{picture}(133.33,123.67)
\put(60.00,10.00){\vector(0,1){110}}
\put(10.00,20.00){\vector(1,0){109.67}}
\linethickness{0.6pt}
\put(40.00,10.00){\line(1,1){78.00}}
\put(10.00,40.00){\line(1,0){60.00}}
\put(70.00,40.00){\line(0,1){80.00}}
\linethickness{1.2pt}
\multiput(55.67,22.33)(0.12,0.29){281}{\line(0,1){0.29}}
\linethickness{0.6pt}
\put(89.20,103.45){\circle*{1.20}}
\put(55.50,22.00){\circle*{1.20}} 
\put(70.00,20.00){\circle*{1.20}}
\put(60.00,40.00){\circle*{1.20}}
\put(57.93,27.93){\circle*{1.20}}
\put(62.90,40.00){\circle*{1.20}}
\put(70.00,57.00){\circle*{1.20}}
\put(70.00,40.00){\circle*{1.20}}
\put(121.33,20.00){\makebox(0,0)[lc]{$\aa(\cdot)$}} 
\put(60.00,121.67){\makebox(0,0)[cb]{$\bb(\cdot)$}} 
\put(72.30,40.00){\makebox(0,0)[lc]{$(\a,\b)$}} 
\put(51.30,23.30){\makebox(0,0)[rc]{$\big(\ax,\bx\xz\big)$}} 
\put(57.14,29.40){\makebox(0,0)[rc]{$(a\zaz{1},b\zaz{1})$}} 
\put(59.47,41.00){\makebox(0,0)[rb]{$(a\zaz{2},b\zaz{2})$}} 
\put(83.53,55.33){\makebox(0,0)[rb]{$(a\zaz{3},b\zaz{3})$}} 
\put(112.67,104.67){\makebox(0,0)[rb]{$\big(\aa(\x'),\bb(\x')\xz\big)$}}
\put(70.00,17.80){\makebox(0,0)[ct]{$\a$}}
\put(59.10,38.80){\makebox(0,0)[rt]{$\b$}}
\put(89.00,35.67){\makebox(0,0)[lt]{$S_1$}}
\put(89.00,76.00){\makebox(0,0)[cb]{$S_2$}}
\put(32.00,31.67){\makebox(0,0)[cb]{$S_3$}}
\put(32.00,76.00){\makebox(0,0)[cb]{$S_4$}}
\put(54.73,22.35){\vector(4,-1){0.2}} %
\multiput(50.97,23.12)(0.27,-0.055){9}{\line(1,0){0.58}}  
\put(62.29,40.5){\vector(3,-1){0.2}} 
\multiput(59.00,41.45)(0.26,-0.075){8}{\line(1,0){0.59}}  
\end{picture}
\caption{An example of line segment
$\big[\big(\ax,\bx,\u01(\x)\xz\big),\big(\aa(\x'),\bb(\x'),$ $\u01(\x')\xz\big)\big]$ in the $\R^3$ space with\xx{ coordinate} axes $\aa(\cdot)$, $\bb(\cdot)$, and $\u01(\cdot)$ projected onto the plane~$\u01=0$.
}
\label{F1}
\end{figure}

Suppose that $(a\zaz{1},b\zaz{1},u\zaz{1}),\ldots,(a\zaz{m},b\zaz{m},u\zaz{m})$, $m\in\{1,2,3\}$, are
the consecutive points where the line segment\xx{ between} $\big[\big(\ax,\bx,\u01(\x)\big),\xy\big(\aa(\x'),\bb(\x'),\u01(\x')\big)\big]$ crosses the\xx{ borders of} planes $\ax=\a,$ $\bx=\b,$ and ${\bx-\ax}=\b-\a$ separating the $S$-regions on the way from $\x$ to~$\x'.$
Then, by the linearity of the segment, it holds that
\eq{a-le}{ 
\ax\le a\zaz{1}\le\cdots\le a\zaz{m}\le \aa(\x'),
}
\eq{b-le}{ 
\bx\le b\zaz{1}\le\cdots\le b\zaz{m}\le \bb(\x'),
}
$$
\u01(\x)<u\zaz{1}<\cdots<u\zaz{m}<\u01(\x')
$$
with strict inequalities in (\ref{a-le}) or in (\ref{b-le}), or in both (since otherwise $\x$ and $\x'$ belong to the same $S$-region).

Consider $\fo$ represented by (\ref{f'''}) as a function $\fhat(a,b,u)$ of $a=\ax,b=\bx$, and $u=\u01(\x)$. Then, using the fact that $\fhat(a,b,u)$ is nondecreasing in all variables on each region, strictly increasing in $u$ on $S_2, S_3,$ and $S_4$, and strictly increasing in $a$ and $b$ on $S_1$ 
and the fact that each point $(a_i,b_i,u_i)$ ($1\le i\le m$) belongs to both regions on the border of which it lies,
we obtain
\begin{eqnarray}\nonumber
\fx
&=&\fhat\big(\ax,\bx,\u01(\x)\big)<\fhat(a\zaz{1},b\zaz{1},u\zaz{1})<\cdots
                             <\fhat(a\zaz{m},b\zaz{m},u\zaz{m})\cr
&<&\fhat\big(\aa(\x'),\bb(\x'),\u01(\x')\big)=f(\x').
\label{fx<fx'} 
\end{eqnarray}

Thus,\xx{ completes the proof of the condition} $\x'\succ\x\To f(\x')>f(\x)$ and\xx{ the proof that} $\fo$ is strictly increasing.\xx{ and the proof of the theorem.} Theorem~\ref{T} is proved.
\qed

\medskip
{\em Proof of Corollary~\ref{c:reduce}.}
$\x\in N\To\Pmx=\Ppx=\emptyset,$ hence $\axl=-\infty$ and $\bxl=+\infty,$ which satisfies the conditions of~$S_4.$

${\axl=\bxl}$ implies $\x\in S_1,$ whence $\fx=\ax$ follows from~(\ref{f'''}).

If $\x\approx\p$ and $\p\in P,$ then since $\fPo$ is gap-safe increasing and thus weakly increasing,
Eq.~(\ref{e:ApprCont}), Remark~\ref{r:=}, and [$(i)\ToTo(iv)\ToTo(v)$] of Proposition~\ref{P:nonstr} imply $\axl=\ap=\fPp=\bp=\bxl,$ hence $\x\in S_1$ and $\fx=\axl=\fPp.$
\qed

\medskip
{\em Proof of Lemma~\ref{l:sin-contour}.}
If $\fP$ is gap-safe increasing, then the conditions presented in Lemma~\ref{l:sin-contour} are satisfied due to Proposition~\ref{P0} and Lemma~\ref{l:strict}.

Conversely, suppose that these conditions\xx{ are satisfied} hold.\xx{ $\fP$ is upper-bounded on {\lowercontour}s and lower-bounded on {\uppercontour}s.}
\xx{Prove that $\fP$ is weakly increasing. }
By the definition of a Pareto set, for any $\p,\,\q\in P,\,$ $\q\succE\p$ reduces to $\q\approx\p,$ and the condition $\big[\xy\q\approx\p\:\To\fPq=\fPp\big]$ implies that $\fP$ is weakly increasing.

Assume that $\fP$ is not gap-safe increasing.
Then there exist $\x,\x'\in\tX$ such that $\x'\succ\x$ and $\bb(\x')\le\ax$. This is possible only if
(a)~$\Pp(\x')=\emptyset$ or
(b)~$\Pmx    =\emptyset$, or
(c)~there are $\p,\,\p'\in P$ such that  $\p'\succE\x'\succ\x\succE\p.$
However,
in (a), $\bb(\x')=+\infty=\ax     $ and $\x \in X$
               (since $\x =\mainfty$ is incompatible with $\ax     =+\infty$ and
                      $\x =\painfty$ is incompatible with $\x'\succ\x$), hence  $\fP$ is not upper-bounded on a \lowercontour.
Similarly,
in (b), $\ax     =-\infty=\bb(\x')$ and $\x'\in X$
               (since $\x'=\painfty$ is incompatible with $\bb(\x')=-\infty$ and
                      $\x'=\mainfty$ is incompatible with $\x'\succ\x$), hence  $\fP$ is not lower-bounded on an \uppercontour.
In (c), by the ``mixed'' strict transitivity of preorders
($\x\succE\y\succ\z\To\x\succ\z$ and $\x\succ\y\succE\z\To\x\succ\z$), we have $\p'\succ\p,$ hence $P$ is not a Pareto set.
In all cases we get a contradiction, therefore, $\fP$ is gap-safe increasing.
\qed

\medskip
{\em Proof of Lemma~\ref{l:setA}.}
Let $\x\in S_1.$ Then $\Ppx\ne\emptyset\;\;\mbox{and\,}\Pmx\ne\emptyset$. Indeed, otherwise either $\bx=+\infty$ or $\ax=-\infty$ and since $\fP$ is upper-bounded on all {\lowercontour}s and lower-bounded on all {\uppercontour}s, $\bx-\ax=+\infty$, which contradicts the assumption. Therefore, $x\in P\cup A.$

Let $x\in P\cup A.$ Then there exist $\p,\q\in P$ such that
\eq{e:2cops}{\q\succE\x\succE\p}
and by the transitivity of $\succE,$ $\q\succE\p.$ Since $P$ is a Pareto set, $\q\not\succN\p.$ 
By the transitivity of $\succN,$ the latter is incompatible with ${\q\succN\x\succN\p}$ in (\ref{e:2cops}), consequently, $\x\approx\p$ for some $\p\in P.$

Let $\x\approx\p$ for some $\p\in P.$ Then by the last statement of Corollary~\ref{c:reduce}, $\x\in S_1.$
This completes the proof.
\qed

\bigskip\bigskip
\appendix
\noindent{\large\bf\appendixname}
\section{Binary relations}
\label{a:BinaryR}

A {\em binary relation\/} $R$ on a set $X$ is a set of ordered pairs $(x,y)$ of elements of $X$ ($R\subseteq X\!\times\! X$);
$(x,y)\in R$ is abbreviated as $xRy$.

A binary relation is
\begin{itemize}
\item
{\em reflexive\/} if $xRx$ holds for every $x\in X$;
\item
{\em irreflexive\/} if $xRx$ holds for no $x\in X$;
\item
{\em transitive\/} if $xRy$ and $yRz$ imply $xRz$ for all $x,y,z\in X$;
\item
{\em symmetric\/} if $xRy$ implies $yRx$ for all $x,y\in X$;
\item
{\em antisymmetric\/} if $xRy$ and $yRx$ imply $x=y$ for all $x,y\in X$;
\item
{\em connected\/} if $xRy$ or $yRx$ holds for all $x,y\in X$ such that $x\ne y.$
\end{itemize}

A binary relation is a/an
\begin{itemize}
\item
{\em preorder\/} if it is transitive and reflexive;
\item
{\em partial order\/} if it is transitive, reflexive and antisymmetric;
\item
{\em strict partial order\/} if it is transitive and irreflexive;
\item
{\em weak order\/} if it is a connected preorder;
\item
{\em linear {\rm(or {\it total\/})} order\/} if it is an antisymmetric weak order (or, equivalently, is a connected partial order);
\item
{\em strict linear order\/} if it is a connected strict partial order;
\item
{\em equivalence relation\/} if it is transitive, reflexive and symmetric.
\end{itemize}

A relation $R$ {\em extends\/} a relation $R_0$ if $R_0\subseteq R.$


\bibliographystyle{splncs04}
\bibliography{all2}

\end{document}